\documentclass[a4paper,12pt]{amsart}
\usepackage{amsmath,amsfonts,amssymb,amsthm,amscd}
\usepackage{latexsym,graphicx}
\usepackage[matrix,arrow,ps,color,line,curve,frame]{xy}
\usepackage[usenames,dvipsnames]{color}

\usepackage[T1]{fontenc}
\usepackage{color}
\usepackage{stmaryrd}

\def\id{{\rm id}}

\renewcommand{\[}{\begin{equation}}
\renewcommand{\]}{\end{equation}}

\newtheorem{thm}{Theorem}[section]

\newtheorem{lem}[thm]{Lemma}
\newtheorem{prop}[thm]{Theorem}

\theoremstyle{definition}

\theoremstyle{remark}

\newcommand{\ra}{\rightarrow}

\newcommand{\nn}{\mathbb{N}}
\newcommand{\cc}{\mathbb{C}}

\newcommand{\h}{\mathfrak{H}}
\newcommand{\ka}{\mathfrak{K}}

\DeclareMathOperator{\Tr}{Tr}

\DeclareMathOperator{\dom}{dom}

\begin{document}
\baselineskip=15pt

\author{Jacopo Bassi}
\address{SISSA (Scuola Internazionale Superiore 
di Studi Avanzati)\\
Via Bonomea 265, 34136 Trieste, Italy
}
\email{jbassi@sissa.it, dabrow@sissa.it}

\author{Ludwik D\k abrowski}

\title[Spectral triples on the Jiang-Su algebra]
{\large Spectral triples on the Jiang-Su algebra}
\begin{abstract}
We construct spectral triples on a class of particular inductive limits of matrix-valued function algebras. 
In the special case  
of the Jiang-Su algebra we employ a particular $AF$-embedding. 
\end{abstract}
\maketitle
\vspace*{-10mm}


\section{Introduction}

According to  the noncommutative differential geometry program  \cite{con, con1} both the topological and the metric
information on a noncommutative space can be  fully encoded as a {\it spectral triple} on the noncommutative algebra of coordinates  on   that space.
Nowadays several noncommutative spectral triples have been constructed, with only a partial
unifying scheme emerging behind some families of examples, e.g. quantum groups and their homogeneous spaces, like quantum spheres and quantum projective spaces.
(see e.g., \cite{ddlw, dd, nt})
Also some preservation properties with respect to the
product, inductive limits or extensions of algebras have been investigated. 

Most of these constructions are still awaiting however a proper analysis of such properties as 
smoothness, dimension (summability)
and other conditions selected by Connes.
As a testing ground for these and related matters
as large as possible class of examples should be investigated, including some important new algebras.

In \cite{qd} 
a general way to construct a spectral triple on arbitrary quasi diagonal $C^*$-algebras was exhibited. 
However, in that case one cannot expect summability. Instead, summability was obtained in  \cite{agi}
for certain inductive family of coverings, and $p$-summability with arbitrary $p$ for any AF-algebra through the construction in \cite{af}.

In the present paper we elaborate a construction that extends the latter mentioned approach to a wider class of particular inductive limits of matrix-valued function algebras whose connecting morphisms have a certain peculiar form.
In particular this construction applies to the Jiang-Su algebra $\mathcal{Z}$ (cf. \cite{js}), 
which was originally constructed 
in terms of an explicit particular inductive limit of dimension drop algebras.
The aim therein was to obtain an example of an infinite-dimensional stably finite nuclear simple unital $C^*$-algebras with exactly one tracial state and with the same $K$-theory of the complex numbers. 
The importance of the Jiang-Su algebra $\mathcal{Z}$ stems from the fact that under some other hypothesis $\mathcal{Z}$-stability entails classification in terms of the Elliott invariant as proved in \cite{winter_2006}.

The organization of the paper is the following: 
In the first section we recall the definition of the Jiang-Su algebra and construct a particular $AF$-embedding for it.
In the second section we compute the image of elements belonging to a dense subalgebra of the Jiang-Su algebra under the representation obtained by composing the forementioned $AF$-embedding with the representation appearing in \cite{af}.
In the last section we use the above results to check that some of the Dirac operators considered in \cite{af} give rise to a spectral triple for the Jiang-Su algebra.

\section{Spectral triple on the Jiang-Su algebra}

Let $B$ be an inductive limit of $C^*$-algebras $B=\lim (B_i ,\phi_i)$, with $B_0 = \cc$ and where every $B_i$ is a unital sub-$C^*$-algebra of the algebra of continuous-valued functions on the interval with values in $M_{n_i}$ for some natural numbers $n_i$ containing a dense $*$-subalgebra of Lipschitz functions and for $l>i$ natural numbers. The connecting morphisms $\phi_{i,i+l}$ take the form

\[
\phi_{i,i+l} (f) = u_{i,i+l} \left( \begin{array}{ccc}	f\circ \xi_{i,i+k}^1 \otimes 1_{N_{i,1}^{i+l}}	&	& 0	\\
																& \ddots	&	\\
															0	&	&	f\circ \xi_{i,ik_i^{i+l}}^{i+l} \otimes 1_{N_{k_i^{i+l}}}	\end{array}	\right)u_{i,i+l}^*
\]

for some natural numbers $k_i^{i+l}$, $N_{i,1}^{i+l}$,...,$N_{i,k_i^{i+l}}$, a unitary $u_{i,i+l} \in C([0,1], M_{n_{i+l}})$ and some paths $\xi_{i,1}^{i+l} , \dots, \xi_{i,k_i^{i+l}}$  satisfying

\begin{equation}
\label{eq}
| \xi_{i,r}^{i+l} (x) - \xi_{i,r}^{i+l} (y) | \leq \frac{1}{2^{l}}, \qquad \mbox{ for } 1\leq r \leq k_i^{i+l}, \; x,y \in [0,1].
\end{equation}

The operators $u_{i,i+l}$ are unitaries in $C([0,1], M_{n_{i+l}})$.\\
The Jiang-Su algebra $\mathcal{Z}$ is an inductive limit of prime dimension drop algebras $Z_i$ satisfying a certain universal property. We will use the original construction appearing in \cite{js}, where it was proven that given $p_i, q_i, n_i= p_i q_i$ defining the prime dimension drop algebra $Z_i$, there are numbers $k_{i,1}^{i+1}$, $k_{i,2}^{i+1}$ and $k_{i,3}^{i+1}$ such that $n_{i+1} =(k_{i,1}^{i+1} + k_{i,2}^{i+1} + k_{i,3}^{i+1})n_i$ is equal to $n_{i+1}=p_{i+1} q_{i+1}$ for some coprime numbers $p_{i+1}$ and $q_{i+1}$ and that there are a unitary $u_{i,i+1} \in C([0,1], M_{n_{i+1}})$ and natural numbers $N_{i,1}^{i+1}$, $N_{i,2}^{i+1}$ and $N_{i,3}^{i+1}$ such that

\[
Z_i \ra Z_{i+1}
\]
\[
\phi_i : f \mapsto u_{i,i+1} \left( \begin{array}{ccc}	f\circ \xi_{i,1}^{i+1} \otimes 1_{N_{i,1}^{i+1}}		&		&	0	\\
												&	f\circ \xi_{i,2}^{i+1} \otimes 1_{N_{i,2}^{i+1}}&		\\
							0					&		&	f\circ \xi_{i,3}^{i+1} \otimes 1_{N_{i,3}^{i+1}}	\end{array}\right)u_{i,i+1}^*
\]

is a connecting morphism for $\xi_1 = x/2$, $\xi_2 = 1/2$ and $\xi_3 = (x+1)/2$.\\
As a consequence, given a natural number $l$, the connecting morphism $Z_i \ra Z_{i+l}$ has the form

\[
\phi_{i,i+l} (f) = u_{i,i+l} \left( \begin{array}{ccc}	f\circ \xi_{i,i+k}^1 \otimes 1_{N_{i,1}^{i+l}}	&	& 0	\\
																& \ddots	&	\\
															0	&	&	f\circ \xi_{i,ik_i^{i+l}}^{i+l} \otimes 1_{N_{k_i^{i+l}}}	\end{array}	\right)u_{i,i+l}^*
\]

for some natural numbers $k_i^{i+l}$, $N_{i,1}^{i+l}$,...,$N_{i,k_i^{i+l}}$, a unitary $u_{i,i+l} \in C([0,1], M_{n_{i+l}})$ and some paths $\xi_{i,1}^{i+l} , ..., \xi_{i,k_i^{i+l}}^{i+l}$ that have the form

\[
\xi_i^{i+l} (x) = 	\frac{x+r}{2^{l}} 	\qquad \mbox{ for } \quad 0\leq r \leq 2^{l} -1 
\]
or
\[
				\xi_i^{i+l} (x) =	\frac{s}{2^{l}}		\qquad \mbox{ for }	\quad 1\leq s \leq 2^{l} -1
					,
\]

It follows that the paths appearing in the connecting morphism $\phi_{l,m}$ satisfy equation \ref{eq} and $\mathcal{Z}$ belongs to the class of inductive limit $C^*$-algebras we want to consider.\\
Note that, given $B$ as above, after reindexing the sequence $B_i$, for example sending $i \mapsto 2i$ we can always suppose that the paths appearing in the connecting morphisms satisfy

\[
|\xi_{i,r}^{i+1}(x) - \xi_{i,r}^{i+1} (y) | \leq \frac{1}{2^i}
\]

for any  $1\leq r \leq k_i^{i+1}$. This relation will be used for the proof of Lemma \ref{iso}.

Fix a sequence of natural numbers $n_i$ as above and consider the inductive limit $A=\lim (A_i , \phi^\circ_i)$, where $A_i = C([0,1], M_{n_i})$ and the connecting morphisms $\phi^\circ_i$ are constructed in the same way as above, but they are considered as unital $*$-homomorphisms between the $A_i$'s. For any $i,l \in \nn$ denote by $\tilde{\phi}^\circ_{i,i+l} : A_i \ra A_{i+l}$ the $*$-homomorphism 

\[
\tilde{\phi}^\circ_{i,i+l} (f) = \left( \begin{array}{ccc}	f\circ \xi_{i,1}^{i+l} \otimes 1_{N_{i,1}^{i+l}}	&	& 0	\\
																& \ddots	&	\\
															0	&	&	f\circ \xi_{i,k_i^{i+l}}^{i+l} \otimes 1_{N_{i,k_i^{i+l}}^{i+l}}	\end{array}	\right).
\]

Let $u_i$ be the unitary corresponding to the connecting morphism $A_1 \ra A_i$ (or $B_1 \ra B_i$).
For any $f \in A_i$ (or $B_i$) there is a unique $\tilde{f} \in A_i$ such that $f=u_i \tilde{f} u_i^*$. In this way the connecting morphisms take the form

\[
\phi^\circ_{i,i+l} (f) = u_{i,i+l} \tilde{\phi}^\circ_{i,i+l}(f) u_{i,i+l}^* = u_{i+l} \tilde{\phi}^\circ_{i,i+l}(\tilde{f}) u_{i+l}^*,
\]

Let now $M= \lim (M_{n_i}, \psi_i)$, where $\psi_i (a) = a \otimes 1_{n_{i+1} /n_i}$.

\begin{lem}
\label{iso}
There is a $*$-isomorphism

\[
\alpha : A \ra M.
\]

Let $\gamma \in (1,2)$. A Lipschitz function $f \in A_i$ with Lipschitz constant $L_f < \gamma^i$ is sent to

\[
\alpha (f) = \lim_{m \ra \infty} \psi_{m}^\infty ( \tilde{\phi}^\circ_{i,m} (\tilde{f}) (0)).
\]

\end{lem}
\proof 
Define $*$-homomorphisms

\[
\alpha_i : A_i \ra M_{n_{i+1}}
\]
\[
f \mapsto \tilde{\phi}^\circ_i(\tilde{f})(0)
\]

and

\[
\beta_i : M_{n_i} \ra A_i
\]
\[
a \mapsto u_{i+1} \bar{a} u_{i+1}^*,
\]
where $\bar{a} \in A_i$ is the constant matrix-valued function taking value $a \in M_{n_i}$. Let now $\gamma \in (1,2)$ and take finite sets $F_i \subset A_i$ consisting of Lipschitz matrix-valued functions with Lipschitz constant less than $\gamma^i$ and such that their union $\bigcup_i F_i$ is dense in $A$. For any $f \in F_i$ and $a \in M_{n_i}$ we have

\[
\alpha_i \circ \beta_i (a) = \psi_i (a),
\]
\[
\| \beta_{i+1} \circ \alpha_i (f) - \phi^\circ_{i,i+1} (f) \| <  \frac{\gamma^i}{2^i}.
\]
Hence the result follows by \cite{classificationnuclear} Proposition 2.3.2. $\Box$\\

\section{The orthogonal decomposition}

Let $\h$ be the Hilbert space considered by Christensen and Antonescu in \cite{af} corresponding to the GNS-representation induced by the unique trace $\tau$ on $M$. This trace is given on the finite-dimensional approximants relative to the inductive limit construction  by the normalized trace on matrices.. Following \cite{af} we want to write $\h$ as an infinite direct sum of the finite dimensional Hilbert spaces on which the $M_{n_i}$'s are represented.\\
Let $\h_i = \overline{M_{n_i}}^\tau$ and let $v \in \h_i$. We can consider $v$ as a matrix of dimension $n_i$ and for any $j <i$, we can write $v$ as a matrix-valued matrix of the form

\[
v=\left( \begin{array}{ccc}	v^{j,i}_{1,1}	&	...		&	v^{j,i}_{1,l_j^i}	\\
					   \vdots		&	\ddots	&	\vdots	\\
					   v^{j,i}_{l_j^i}	&	...		&	v^{j,i}_{l_j^i , l_j^i}	\end{array}\right),
\]

where $l_j^i= n_i / n_j$ is the multiplicity of $M_{n_j}$ in $M_{n_i}$ and the $v^{j,i}_{k,l}$ are matrices in $M_{n_j}$; in particular we can apply the same procedure to these matrices by iteration. With this notation, the projection $P_{i,j}$ from $\h_i$ to $\h_j$ reads

\[
P_{i,j} (v) = \frac{1}{l_j^i}\sum_{k=1}^{l_j^i} v^{j,i}_{k,k} \quad \in M_{n_j}.
\]

If $i>1$, the projection $R_j$ from $\h_j$ to the orthogonal complement of $\h_{j-1}$ in $\h_j$ reads for $w \in \h_j$:

\[
R_j (w) = \]
\[
\left( \begin{array}{cccc}	w^{j-1,j}_{1,1}-\frac{1}{l_{j-1}^j}\sum_{k=1}^{l_{i-1}^i}w^{j-1,j}_{k,k}	&	w^{j-1,j}_{1,2}	&	...	&	w^{j-1,j}_{1,l_{j-1}^j}	\\
							w^{j-1,j}_{2,1}		&	w^{j-1,j}_{2,2} - \frac{1}{l_{j-1}^j}\sum_{k=1}^{l_{i-1}^i}w^{j-1,j}_{k,k}	&			&	w^{j-1,j}_{2,l_{j-1}^j}	\\
							\vdots	&	&	\ddots	&		\vdots	\\
							w^{j-1,j}_{l_{j-1}^j,1}	&	...	&	&	w^{j-1,j}_{l_{j-1}^j , l_{j-1}^j} -\frac{1}{l_{j-1}^j}\sum_{k=1}^{l_{i-1}^i}w^{j-1,j}_{k,k}	\end{array}\right).
\]

Hence, if we denote by $\ka_i = \h_i \ominus \h_{i-1}$, the projection $Q_j : \h \ra \mathfrak{K}_j$, when applied to an element $v \in \h_i$ takes the form, for $1\leq s,t \leq l_{j-1}^j$

\[
(Q_j (v))_{s,t}^{j-1,j} = \begin{cases}	\frac{1}{l_{j}^i}\sum_{k=1}^{l_j^i} (v^{j,i}_{k,k})^{j-1,j}_{s,s} -\frac{1}{l_{j-1}^i}\sum_{t=1}^{l_{j-1}^j}\sum_{k=1}^{l_j^i} (v^{j,i}_{k,k})^{j-1,j}_{t,t} 	&	\mbox{ for } s=t \\
	\frac{1}{l_j^i}\sum_{k=1}^{l_j^i} (v^{j,i}_{k,k})^{j-1,j}_{s,t}	&	\mbox{ for } s\neq t
	\end{cases}.
	\]

\section{The commutators}

Take $i <n<m$ and $v \in \h_m$, $f \in A_i$. We want to compute the elements $Q_n (\tilde{\phi}^\circ_{i,m} (\tilde{f}) (0) v)$ and $\tilde{\phi}^\circ_{i,n} (\tilde{f})(0) Q_n v$. To this end we want to write $\tilde{\phi}^\circ_{i,m} (\tilde{f})$ as the composition $\tilde{\phi}^\circ_{n,m} \circ \tilde{\phi}^\circ_{n-1,n} \circ \tilde{\phi}^\circ_{i,n-1} (\tilde{f})$.\\
Let $k_j^i$ be the amount of different paths appearing in the connecting morphism $\phi_{j,i}$. If $1\leq j \leq k_{n-1}^n$, we denote by $\tilde{f}  \circ [\xi_i^{n-1}] \circ \xi_{n-1,j}^n = \tilde{\phi}^\circ_{i,n-1}(\tilde{f}) \circ \xi_{n-1,j}^n$ the matrix-valued function

\[
\left( \begin{array}{ccc} \tilde{f}  \circ  \xi_{i,1}^{n-1} \circ \xi_{n-1 , j}^n\otimes 1_{N_{i,1}^{n-1}} 	&	&	0\\
									&	\ddots		&	\\
								0	&				&	\tilde{f}   \circ \xi_{i,k_i^{n-1}}^{n-1}\circ \xi_{n-1,j}^n\otimes 1_{N_{i,k_i^{n-1}}^{n-1}}	\end{array} \right),
\]

then we can write 

\[
\tilde{\phi}^{\circ}_{i,n} (\tilde{f}) = \tilde{\phi}^{\circ}_{n,n-1} \circ \tilde{\phi}^{\circ}_{i,n-1} (\tilde{f}) =
\]
\[
\left(\begin{array}{ccc}	\tilde{f}  \circ [\xi_i^{n-1}] \circ \xi_{n-1 , 1}^n \otimes 1_{N_{n-1 , 1}^n }	&	&	0	\\
												&	\ddots	&	\\
										0	&	&	\tilde{f}  \circ [\xi_i^{n-1}] \circ \xi_{n-1 , k_{n-1}^n}^n \otimes 1_{N_{n-1 , k_{n-1}^n}^n} \end{array}\right).
\]

For $1\leq s, \leq l_{n-1}^n$, we denote by $\bar{\xi}_{n-1,s}^n$ the path

\[
\bar{\xi}_{n-1,s}^n = \begin{cases}	\xi_{n-1,1}^n 	&\mbox{ for } 1\leq s \leq N_{n-1,1}^n\\
							\xi_{n-1,2}^n	&\mbox{ for } N_{n-1,1}^n < s \leq N_{n-1, 1}^n + N_{n-1,2}^n \\
							\vdots  & \\
							\xi_{n-1,k_{n-1}^n} 	& \mbox{ for } \sum_{k=1}^{k_{n-1}^n - 1}N_{n-1 , k}^n< s \leq l_{n-1}^n
							\end{cases}.
\]

Thus we obtain for $1\leq s,t \leq l_{n-1}^n$,

\[
(\tilde{\phi}^\circ_{n-1,n} \circ \tilde{\phi}^\circ_{i,n-1} (\tilde{f})(0) Q_n v)^{n-1,n}_{s,t} =
\]
\[
\frac{1}{l_n^m} \sum_{j=1}^{l_n^m} (\tilde{f}  \circ [\xi_i^{n-1}] \circ \bar{\xi}_{n-1,s}^n)(0) (v_{j,j}^{n,m} )_{s,t}^{n-1,n} \qquad \mbox{ for } s\neq t
\]
and 
\[
\frac{1}{l_n^m} \sum_{j=1}^{l_n^m} (f \circ [\xi_i^{n-1} ]\circ \bar{\xi}_{n-1,s}^n )(0) \left( (v_{j,j}^{n,m})_{s,s}^{n-1,n} - \frac{1}{l_{n-1}^n} \sum_{k=1}^{l_{n-1}^n} (v_{j,j}^{n,m})_{k,k}^{n-1,n}\right) 	\quad \mbox{ for } s=t.
\]

In the same way, for $1\leq j \leq l_n^m$, we can define paths

\[
\bar{\xi}_{n,j} = \begin{cases}	\xi_{n,1}^m	&\mbox{ for }	1\leq j \leq N_{n,1}^m	\\
						\xi_{n,2}^m	&\mbox{ for } 	N_{n,1}^m < j \leq N_{n,1}^m + N_{n,2}^m\\
						\vdots	&	\\
						\xi_{n,k_n^m}^m 	&\mbox{ for }	\sum_{k=1}^{k_n^m -1} N_{n,k}^m < j \leq l_n^m
						\end{cases}
						\]

and compute for $1\leq s,t \leq l_{n-1}^n$,

\[
(Q_n \tilde{\phi}^{\circ}_{i,m} (\tilde{f}) (0) v)_{s,t}^{n-1,n}= (Q_n (\tilde{\phi}^{\circ}_{n,m}  \circ \tilde{\phi}^{\circ}_{n-1,n}  \circ \tilde{\phi}^{\circ}_{i,n-1})(\tilde{f}) (0) v)_{s,t}^{n-1,n}=
\]
\[
\frac{1}{l_n^m} \sum_{j=1}^{l_n^m} (\tilde{f}   \circ [ \xi_i^{n-1}]\circ \bar{\xi}_{n-1,s}^n\circ \bar{\xi}_{n,j}^m) (0) (v_{j,j}^{n,m})_{s,t}^{n-1,n} \quad \mbox{ for } s \neq t
\]

and 

\[
\frac{1}{l_n^m} \sum_{j=1}^{l_n^m} [ ( \tilde{f}  \circ [\xi_i^{n-1}]\circ \bar{\xi}_{n-1 , s}^n \circ \bar{\xi}_{n,j}^m)(0) (v_{j,j}^{n,m} )_{s,s}^{n-1 , n}
\]
\[
- \frac{1}{l_{n-1}^n} \sum_{k=1}^{l_{n-1}^n} (\tilde{f}   \circ [\xi_i^{n-1} ]\circ \bar{\xi}_{n-1,k}^n\circ \bar{\xi}_{n,j}^m)(0) (v_{j,j}^{n,m} )_{k,k}^{n-1 , n} ) ]\quad \mbox{ for } s=t.
\]

Thus we can write the commutators

\[
(Q_n (\tilde{\phi}^{\circ}_{i,m} (\tilde{f}) (0) v)- \tilde{\phi}^{\circ}_{i,n} (\tilde{f})(0) Q_n v)_{s,t}^{n-1,n}=
\]
\[
\frac{1}{l_n^m} \sum_{j=1}^{l_n^m} (\tilde{f}   \circ [\xi_i^{n-1}]\circ \bar{\xi}_{n-1,s}^n \circ \bar{\xi}_{n,j}^m - \tilde{f}  \circ [\xi_i^{n-1}]\circ \bar{\xi}_{n-1,s}^n )(0)(v_{j,j}^{n,m})_{s,t}^{n-1,n} \quad \mbox{ for } s\neq t
\]
and
\[
\frac{1}{l_n^m} \sum_{j=1}^{l_n^m} [ (\tilde{f} \circ [\xi_i^{n-1}] \circ \bar{\xi}_{n-1,s}^n \circ \bar{\xi}_{n,j}^m - \tilde{f} \circ [\xi_i^{n-1}]  \circ \bar{\xi}_{n-1,s}^n )(0) (v_{j,j}^{n,m})_{s,t}^{n-1,n} +\]
\[
\frac{1}{l_{n-1}^n} \sum_{k=1}^{l_{n-1}^n} (\tilde{f}  \circ [\xi_i^{n-1}]\circ \bar{\xi}_{n-1,s}^n - \tilde{f}   \circ [\xi_i^{n-1}]\circ \bar{\xi}_{n-1,k}^n \circ \xi_{n,j}^m )(0) (v_{j,j}^{n,m})_{k,k}^{n-1,n} \circ \bar{\xi}_{n,j}^m] \quad \mbox{ for } s=t.
\]

\begin{lem}
\label{lip}
Let $i<l<m\leq k$ be natural numbers and let $\xi_i^l$, $\xi_l^m$, $\xi_l^k$ be paths on the interval $[0,1]$ such that

\[
| \xi_{i}^{l} (x) - \xi_{i}^{l} (y) | \leq \frac{1}{2^{l-i}}, \qquad \mbox{ for any } x,y \in [0,1].
\]

Then, given any $n >0$ and any Lipschitz function in $C([0,1], M_n)$ with Lipschitz constant $L_f$, we have

\[
\| (f \circ  \xi_i^l \circ  \xi_l^m)(0) - (f\circ \xi_i^l   \circ \xi_l^k)(0)\| \leq \frac{2^{i} L_f}{2^l}.
\]
\end{lem}
\proof 

This is a consequence of the fact that $|\xi_i^l (x) - \xi_i^l (y)| \leq \frac{1}{2^{l-i}}$ for every $x,y \in [0,1]$. $\Box$\\

\section{The spectral triple}
Note that if $D= \sum_n \alpha_n Q_n$ for a certain sequence of real numbers 
$\{ \alpha_n \}$, then the domain of $D$, 
$\dom (D) = \{ v \in \h \; : \; \{ \|\alpha_n Q_n v\|\} \in l^2 (\nn)\}$ 
is left invariant under the action of any $f \in A$; thus in particular, for every $f \in B$ and it makes sense to consider the (in general unbounded) operator $[D,f]$.\\
Moreover, it follows from the Hann-Banach extension Theorem, that if $T$ is an unbounded operator on $\h$ whose domain contains the algebraic direct sum $\oplus_{alg} \ka_i$ and $\| TP_n\|$ is uniformly bounded on $n$, then $T$ extends (uniquely) to a bounded operator on the whole Hilbert space $\h$.\\
Hence, to obtain boundedness of $[D,f]$, we want to compute estimates for $\|[D,f]P_n\|$ for every $n$.\\
For every $i \in \nn$ we will denote by $LB_i$ the linear subspace of $B_i$ consisting of Lipschitz functions with Lipschitz constant smaller than $\gamma^{i}$ for some $\gamma \in (1,2)$. Observe that $\phi^\circ |_{LB_i}$ is a linear map sending $LB_i$ into $LB_{i+1}$ and that the algebraic direct limit $\bigcup_i LB_i$ is a dense $*$-subalgebra of $B$.

\begin{prop}
Let $D=\sum_n \alpha_n Q_n$, with $\{ \alpha_n \}$ a diverging sequence of real numbers satisfying $\alpha_0 = 0$, $|\alpha_n | \leq \beta^{2(n-1)}$ with $\beta <2$ and $n >0$. Then $(\bigcup_i LB_i , \h , D)$ is a spectral triple for $B$.\\
It is is $p$-summable whenever the sequences of numbers $\{ \alpha_i \}$, $\{ n_i\}$ satisfy
\[
 \sum_{i \geq 1} (1+ \alpha_i^2)^{-p/2} (n_i^2 - n_{i-1}^2) <\infty
\]
for some $p >0$.
\end{prop}
\proof 
After reindexing $i \mapsto 2i$, the $*$-isomorphism $\alpha: A \ra M$ has the concrete description given in Lemma \ref{iso}. Thus we can compose it with the GNS representation of $M$ induced by the unique trace $\tau$.\\
Let $l,m \in \nn$ and $v \in \h_l$. Denote by $\beta_{l,m}^\h : \h_l \ra \h_m$ and $\beta_{l,\infty}^\h : \h_l \ra \h$ the connecting isometries. Note that for $i<n \in \nn$ and $f \in LB_i$ the action of $f$ on $v$ reads

\[
\lim_{m \ra \infty} \beta_{m,\infty}^\h \tilde{\phi}_{i,m}^\circ (\tilde{f})(0) \beta_{l,m}^\h v,
\]

where we use the convention that $\tilde{\phi}_{i,m}^\circ = \id$ for $m \leq i$ and $\beta_{l,m}^\h = \id$ for $m \leq l$.
Thus we can write

\[
\| Q_n f v - fQ_n v\| = \| \beta_n^{\h,\infty} Q_n \lim_{m\ra \infty} \tilde{\phi}^\circ_{i,m} (\tilde{f}) (0) \beta_{l,m}^\h v - \lim_{m\ra \infty} \beta_{m,\infty}^\h \tilde{\phi}^\circ_{i,m} (\tilde{f})(0) \beta_{n,m}^\h Q_n v\|.
\]

Since the sequence $\beta_{m,\infty}^\h \tilde{\phi}^\circ_{i,m} (f)(0) \beta_{l,m}^\h v$ converges, there is an $M$ such that 

\[
\| \beta_{k,\infty}^\h \tilde{\phi}^\circ_{i,k} (\tilde{f})(0) \beta_{l,k}^\h v - \lim_{m\ra \infty} \beta_{m,\infty}^\h \tilde{\phi}^\circ_{i,m} (\tilde{f}) (0) \beta_{l,m}^\h v\| \leq \frac{1}{2^{2(n-1)}}
\]

for any $k \geq M$.
Moreover, by Lemma \ref{lip} and the discussion preceding it

\[
\begin{split}
\| [\beta_{n,m}^\h \tilde{\phi}^\circ_{i,n}(\tilde{f}) (0) -& \tilde{\phi}^\circ_{i,m} (f)(0) \beta_{n,m}^\h]Q_n v\|\\
& =\| (\beta_{n,m}^\h \tilde{f}\circ [\xi_i^n] (0) -  \tilde{f} \circ [\xi_i^n] \circ [\xi_n^m] (0)\beta_{n,m}^\h)Q_n v\| \leq \frac{2^{2i} L_f}{2^{2(n-1)}}
\end{split}
\]

for $m >n$ and
\[
 \| Q_n \tilde{\phi}^\circ_{i,M}(\tilde{f})(0) \beta_{l,M}^\h v - \tilde{\phi}^\circ_{i,n} (\tilde{f})(0) Q_n \beta_{l,M}^\h v\| \leq \frac{2^{2i} L_f}{2^{2(n-1)}} .
\]
We can suppose $M >n$ and obtain 

\[
\begin{split}
\| \beta_{n,\infty}^\h Q_n   \lim_{m\ra \infty} \beta_{m,\infty}^\h\tilde{\phi}^\circ_{i,m}& (\tilde{f}) (0) \beta_{l,m}^\h v -  \lim_{m\ra \infty} \beta_{m,\infty}^\h \tilde{\phi}^\circ_{i,m} (\tilde{f})(0) \beta_{n,m}^\h Q_n v\| \\
&\leq \|\beta_{n,\infty}^\h Q_n [\beta_{M,\infty}^\h\tilde{\phi}^\circ_{i,M} (\tilde{f})(0) \beta_{l,M}^\h v - \lim_{m\ra \infty}\beta_{m,\infty}^\h \tilde{\phi}^\circ_{i,m} (\tilde{f}) (0) \beta_{l,m}^\h v]\| \\
& + \| Q_n \tilde{\phi}^\circ_{i,M}(\tilde{f})(0) \beta_{l,M}^\h v - \tilde{\phi}^\circ_{i,n} (\tilde{f})(0) Q_n \beta_{l,M}^\h v\| \\
& + \| \lim_{m \ra \infty} \beta_{m,\infty}^\h [ \beta_{n,m}^\h \tilde{\phi}^\circ_{i,n} (\tilde{f}) (0) - \tilde{\phi}^\circ_{i,m} (\tilde{f})(0) \beta_{n,m}^\h ] Q_n \beta_{l,M}^\h v \| \\
&\leq \frac{1 + 2^{2i+1} L_f}{2^{2(n-1)}}.
\end{split}
\]
Thus we obtain
\[
\| [\alpha_n Q_n,f] P_m\| \leq \frac{|\alpha_n |( 1+2^{2i+1} L_f)}{2^{2(n-1)}} \leq  (1+2^{2i+1} L_f)(\beta /2)^{2(n-1)}.
\]

Hence 
\[
\begin{split}
\|[D,f]\|& \leq \| [\sum_{n=1}^i \alpha_n Q_n ,f \| + \| \sum_{n>i} \alpha_n Q_n , f ]\| \\
				&\leq 2 \|f\| \sum_{n=1}^i |\alpha_n | + ( 1+2^{2i+1} L_f)\sum_{n >i}  (\beta/2)^{2(n-1)} < \infty
\end{split}
\]
and $[D,f]$ extends to a bounded operator.\\
Moreover $D$ has compact resolvent since it has discrete spectrum and its eigenvalues have finite multiplicity.
Suppose we have sequences $\{\alpha_i\}$, $\{n_i\}$ and a real number $p>0$ as in the statement. Then

\[\label{dim}
\Tr ( (1+D^2)^{-{p/2}}) = 1 + \sum_{i\geq 1} (1 + \alpha_i^2)^{-p/2} (n_i^2 - n_{i-1}^2) < \infty. \qquad \Box
\]

As the final comment we observe that by looking at the growth of the dimensions of the matrix algebras appearing in the original construction of the Jiang-Su algebra (cfr.  \cite{js}), it is clear that \eqref{dim} can not be satisfied 
and the spectral triples exhibited above are not $p$-summable. 
Also, with the help of Stirling formula it can be seen that Tr\,$\exp{(-D^2)}$ diverges and thus the $\theta$-summability does not hold either.

\section*{Acknowledgments}
\noindent
This work is part of the project 
\emph{Quantum Dynamics} sponsored by  EU-grant RISE~691246. J.B. thanks Prof. Wilhelm Winter for the hospitality at the University of Munster.
L.D. is grateful for the support at IMPAN provided by
Simons-Foundation grant 346300 and a Polish Government MNiSW 2015--2019 matching fund.


\end{document}